\documentclass[12pt,a4paper]{article}

\usepackage{amsmath}
\usepackage{amssymb}
\usepackage{color}
\usepackage[dvips]{graphics}
\usepackage{cite}

\newtheorem{theorem}{Theorem}
\newtheorem{lemma}{Lemma}
\newtheorem{corol}{Corollary}

\newcommand{\fin}{\,\rule{1ex}{1.6ex}\,}

\date{~}
\definecolor{light}{gray}{.85}

\title{Multistep $\varepsilon$--algorithm,
Shanks' transformation, and Lotka--Volterra system by Hirota's method}

\author{Claude Brezinski\thanks{Laboratoire Paul Painlev\'e, UMR CNRS 8524, UFR
de Math\'ematiques Pures et Appliqu\'ees, Universit\'e des Sciences
et Technologies de Lille,
France, E--mail: {\tt Claude.Brezinski@univ-lille1.fr}.}
\and Yi He\thanks{LSEC, Institute of Computational Mathematics
and Scientific Engineering Computing, AMSS, Chinese Academy of
Sciences,
and
Graduate School of the Chinese Academy of Sciences, Beijing, PR China. E--mail: {\tt \{heyi, sunjq\}@lsec.cc.ac.cn},}
\and Xing-Biao Hu\thanks{LSEC, Institute of Computational Mathematics
and Scientific Engineering Computing, AMSS, Chinese Academy of
Sciences,
Beijing, PR China, E--mail: {\tt hxb@lsec.cc.ac.cn}.}
\and Michela Redivo--Zaglia\thanks{Universit\`a degli Studi di Padova,
Dipartimento di Matematica Pura ed Applicata,
Italy. E--mail : {\tt Michela.RedivoZaglia@unipd.it}.}
\and Jian-Qing Sun$^\dag$
}

\begin{document}

\maketitle

\begin{abstract}
In this paper, we give a multistep extension of the $\varepsilon$--algorithm of Wynn, and we show that it implements
a multistep extension of the Shanks' sequence transformation which is defined by ratios of determinants.
Reciprocally, the quantities defined in this transformation can be recursively computed by the multistep
$\varepsilon$--algorithm. The multistep $\varepsilon$--algorithm and the multistep Shanks' transformation
are related to an extended discrete Lotka--Volterra system.
These results are obtained by using the Hirota's bilinear method, a procedure
quite useful in the solution of nonlinear partial differential and difference equations.
\end{abstract}

\section{The scenery}
\label{sce}

Let $(S_n)$ be a sequence of numbers converging to $S$. If its convergence is slow, it can be transformed, by
a {\it sequence transformation}, into a set of new sequences $\{(T_k^{(n)})\}$, depending on two indexes $k$
and $n$, and converging, under certain assumptions, faster to the same limit, that is such that
$$\lim_{n \to \infty} \frac{T_k^{(n)}-S}{S_n-S}=0, \quad \mbox{or} \quad \lim_{k \to \infty} \frac{T_k^{(n)}-S}{S_k-S}=0,
\quad \mbox{or both}.$$

A well--known example of such a transformation is the Richardson extrapolation process, which gives rise to
the Romberg's method for accelerating the convergence of the trapezoidal rule for approximating a definite integral.
Let us mention that sequence transformations can also be applied to diverging power sequences,
thus leading, in some situations, to interesting results such as analytic continuation (this is the case of
the $\varepsilon$--algorithm which,
applied to the partial sum of a divergent power series, computes its Pad\'{e} approximants).

\vskip 2mm

In many sequence transformations, the terms of the new sequences can be expressed as ratios of determinants,
and there exists, in each particular case, a (usually nonlinear) {\it recursive algorithm} for avoiding the computation
of these determinants and implementing the transformation under consideration \cite{cbmrz,sidi,weni,walz,wimp}.

\vskip 2mm

The most well--known transformation of this type is due to Shanks \cite{ds,shanks}. It
can be implemented via the $\varepsilon$--algorithm of Wynn \cite{wynn}.
Recently, a new recursive algorithm for accelerating the convergence of sequences was derived by
He, Hu, Sun and Weniger \cite{huweni} from the lattice Boussinesq equation. This algorithm resembles to
the $\varepsilon$--algorithm, and it was proved that the quantities it computes can be expressed as ratios
of determinants, thus extending the Shanks' sequence transformation.
In this paper, inspired by this approach, we will extend further the $\varepsilon$--algorithm,
and we will show that it implements an extension of the Shanks' transformation,
thus leading to a multistep $\varepsilon$--algorithm and a multistep Shanks' transformation.
The proof makes use of the Hirota's bilinear method \cite{hiro} which was invented for resolving integrable
nonlinear partial differential or difference evolution equations having soliton solutions.

\vskip 2mm

For some years now, there has been a great concern for convergence acceleration algorithms among the
community of mathematical physicists working on integrable systems, KdV and other equations, soliton theory,
Toda lattices, etc. \cite{cr,naga,nagat,int4,int3}. These researchers are interested by the fact that
convergence acceleration algorithms are
nonlinear difference equations in two variables whose solutions are explicitly known. Determinants often play a
central role in this type of problems as exemplified, for example,  in \cite{vein}.
An important procedure for obtaining a closed--form solution of soliton equations is the Hirota's bilinear method \cite{hiro}
which consists in writing the solution as a ratio, and then working with its numerator and its denominator.

\vskip 2mm

In Section \ref{sshanks}, we discuss the Shanks' sequence transformation and its implementation by
the $\varepsilon$--algorithm of Wynn. The quantities involved in this transformation and in this algorithm
are expressed by ratios of Hankel determinants. In Section \ref{meps}, we present our multistep extension of the
$\varepsilon$--algorithm, and the corresponding multistep extension of the Shanks' transformation.
Section \ref{sshan} is devoted to some relations between determinants that will be useful for our purpose.
The Hirota's bilinear method is presented in Section \ref{sshiro}. In Section \ref{sse}, we show that the quantities recursively
computed by the multistep $\varepsilon$--algorithm correspond to the ratios of determinants defining the multistep
Shanks' transformation, and, reciprocally, in Section \ref{sss}, we show that the multistep Shanks' transformation can be
implemented by the multistep $\varepsilon$--algorithm. Finally, in Section \ref{slv}, the connection between an extended
discrete hungry Lotka--Volterra system and the multistep $\varepsilon$--algorithm is discussed.
Hirota's method is essential for obtaining these results.
The paper ends by some considerations on further researches.

\section{The Shanks' transformation and the $\varepsilon$--algorithm}
\label{sshanks}

The Shanks' sequence transformation \cite{ds,shanks} $e_k : (S_n)  \longmapsto \{(e_k(S_n))\}$ consists in transforming
a given sequence $(S_n)$ into the set of sequences $\{(e_k(S_n))\}$ whose terms are defined by
\begin{equation}
e_k(S_n)=\frac{{\cal H}_{k+1}(S_n)}{{\cal H}_k(\Delta^2 S_n)}, \quad k,n=0,1,\ldots, \label{shanks}
\end{equation}
where $\Delta$ is the usual forward difference operator whose powers are defined by
$$\Delta^{i+1} S_n=\Delta^i S_{n+1}-\Delta^i S_n$$
with $\Delta^0 S_n=S_n$, and where ${\cal H}_k(u_n)$ denotes the Hankel determinant
\begin{equation}
{\cal H}_k(u_n)=
\left|
\begin{array}{cccc}
u_n &u_ {n+1} & \cdots & u_{n+k-1}\\
u_{n+1} &u_ {n+2} & \cdots & u_{n+k}\\
\vdots & \vdots && \vdots \\
u_{n+k-1} &u_ {n+k} & \cdots & u_{n+2k-2}
\end{array}
\right|, \nonumber
\end{equation}
with ${\cal H}_0(u_n)=1$.

Obviously, replacing each row, in this determinant,  by its difference with the previous one, repeating this
operation several times, and performing it also on the columns, we have
\begin{equation}
{\cal H}_k(u_n)=
\left|
\begin{array}{ccc}
u_n & \cdots & u_{n+k-1}\\
\Delta u_{n}  & \cdots & \Delta u_{n+k-1}\\
\Delta^2 u_{n}  & \cdots & \Delta^2 u_{n+k-1}\\
\vdots  && \vdots \\
\Delta^{k-1} u_{n} & \cdots & \Delta^{k-1} u_{n+k-1}
\end{array}
\right|=
\left|
\begin{array}{cccc}
u_n & \Delta u_{n} &\cdots & \Delta^{k-1} u_{n} \\
\Delta u_{n}  & \Delta^2 u_{n} & \cdots & \Delta^k u_{n}\\
\vdots && \vdots \\
\Delta^{k-1} u_{n} & \Delta^{k} u_{n} & \cdots & \Delta^{2k-2} u_{n}
\end{array}
\right|. \nonumber
\end{equation}

The $\varepsilon$--algorithm is a recursive algorithm due to Wynn \cite{wynn} for implementing the Shanks' transformation
without computing the Hankel determinants appearing in (\ref{shanks}). Its rule is
\begin{equation}
\varepsilon_{k+1}^{(n)}=\varepsilon_{k-1}^{(n+1)}+\frac{1}{\varepsilon_{k}^{(n+1)}-\varepsilon_{k}^{(n)}},
 \qquad k,n=0,1,\ldots
\label{epsi}
\end{equation}
with $\varepsilon_{-1}^{(n)}=0$ and $\varepsilon_{0}^{(n)}=S_n$, for $n=0,1,\ldots$

The connection between the $\varepsilon$--algorithm and the Shanks' transformation is given by
\begin{equation}
\varepsilon_{2k}^{(n)}=e_k(S_n) \quad \mbox{and} \quad \varepsilon_{2k+1}^{(n)}=\frac{1}{e_k(\Delta S_n)},
\quad k,n=0,1,\ldots \label{se}
\end{equation}
Thus, the $\varepsilon_{2k+1}^{(n)}$'s are intermediate results, and we have
\begin{equation}
\varepsilon_{2k}^{(n)}=\frac{{\cal H}_{k+1}(S_n)}{{\cal H}_k(\Delta^2 S_n)} \quad \mbox{and} \quad
\varepsilon_{2k+1}^{(n)}=\frac{{\cal H}_k(\Delta^3 S_n)}{{\cal H}_{k+1}(\Delta S_n)}. \nonumber
\end{equation}

\vskip 2mm

The quantities $\varepsilon_{k}^{(n)}$ are usually displayed in a two--dimensional array (the $\varepsilon${\it --array})
where the lower index $k$ remains the
same in a column of the table, and the upper index $n$ is the same in a descending diagonal. Thus, the rule
(\ref{epsi}) relates four quantities located at the four vertices of a lozenge in three different columns and two
descending diagonals as showed below

$$
\begin{array}{cccc}
&\varepsilon_{k}^{(n)}&&\\
\varepsilon_{k-1}^{(n+1)} && \varepsilon_{k+1}^{(n)}&\\
&\varepsilon_{k}^{(n+1)}&&
\end{array}
$$

For implementing the $\varepsilon$--algorithm efficiently, the best technique, due to Wynn \cite{asc1,asc2},
consists in storing the last ascending
diagonal of the $\varepsilon$--array (in this diagonal the sum of the lower and the upper indexes is constant),
and to add, one by one, the terms of the sequence to be transformed. Then, a new ascending diagonal is built
step--by--step, by moving up the lozenge, and the new diagonal gradually replaces the old one. The
corresponding {\sc fortran} subroutine can be found in \cite{cbmrz}.

\vskip 2mm

Since the quantities with an odd lower index are intermediate computations, they can be eliminated, thus leading
to the {\it cross rule} also due to Wynn \cite{cross}
$$\frac{1}{\varepsilon_{2k+4}^{(n)}-\varepsilon_{2k+2}^{(n+1)}}+
\frac{1}{\varepsilon_{2k}^{(n+2)}-\varepsilon_{2k+2}^{(n+1)}}
=\frac{1}{\varepsilon_{2k+2}^{(n+2)}-\varepsilon_{2k+2}^{(n+1)}}+
\frac{1}{\varepsilon_{2k+2}^{(n)}-\varepsilon_{2k+2}^{(n+1)}},$$
with the initial conditions $\varepsilon_{-2}^{(n)}=\infty$ and $\varepsilon_0^{(n)}=S_n$ for $n=0,1,\ldots$
Obviously, it is also possible to eliminate the $\varepsilon_k^{(n)}$'s with an even lower index for obtaining a rule only
involving quantities with a lower odd index, although this is less useful from the numerical point of view.

\vskip 2mm

The proof given by Wynn for his $\varepsilon$--algorithm was mostly a {\it verification} of the link between
the Shanks' transformation and the algorithm, since he introduced the ratios of Hankel determinants for $e_k(S_n)$ and
$e_k(\Delta S_n)$ into the rule of the $\varepsilon$--algorithm, and he showed that the equality held by making use
of the Sylvester's determinantal identity and the Schweins' one which can be found, for example, in \cite{aitk} (see
\cite[pp. 142--143]{bior} for their proofs). The difficulty of the proof resided in the nonlinearity of the algorithm.
Of course, Wynn's great merit was the idea of the $\varepsilon$--algorithm itself, followed by this verification.

\vskip 2mm

There are three approaches for linking a sequence transformation and a (usually nonlinear) recursive algorithm for
its implementation. By increasing order of complexity, they are
\begin{enumerate}
\item {\it Verification:} the transformation and the algorithm are both known, and one has to {\it verify} that they
lead to identical sequences. This is the way followed by Wynn in \cite{wynn} when he gave his $\varepsilon$--algorithm.
\item {\it Derivation:} only the transformation is known, and one has to {\it derive} an algorithm for its implementation.
This is the case,
for example, of the $E$--transformation which is the most general sequence transformation known
so far, and which can be implemented by the $E$--algorithm, an algorithm which appeared
almost simultaneously  in various contexts \cite{ealg,5havi,5mein,5schn}.
This was also certainly the way  Wynn followed when he derived his $\varepsilon$--algorithm, although it was not
presented like that in his paper \cite{wynn}.
\item {\it Proof:} only the algorithm is known, and one has to guess a formula (that is a ratio of determinants)
for the transformation it is implementing, and to {\it prove} it.
This was the situation for the second generalization of the $\varepsilon$--algorithm proposed in \cite{geneps},
whose form was obtained by Salam \cite{salam,salam2}.
Let us mention that the $\theta$--algorithm \cite{theta} is an extrapolation algorithm for which no determinantal
formula is known yet, if it exists.
\end{enumerate}

Now, after presenting the multistep $\varepsilon$--algorithm and the multistep Shanks' transformation
(Section \ref{meps}), we will
show, with the help of determinantal identities (Section \ref{sshan}) and the Hirota's bilinear method
(Section \ref{sshiro}), how to go from the multistep
$\varepsilon$--algorithm to the multistep Shanks' transformation  (Section \ref{sse}), and back  (Section \ref{sss}).

\section{The multistep $\varepsilon$--algorithm and the multistep Shanks' transformation}
\label{meps}

Let $m$ be a fixed strictly positive integer. We define the {\it multistep $\varepsilon$--algorithm} by the recursive rule
\begin{equation}
\label{gd}
\varepsilon_{k+1,m}^{(n)}=\varepsilon_{k-m,m}^{(n+1)}+\frac{1}{\prod_{i=1}^{m}(\varepsilon_{k-m+i,m}^{(n+1)}-
\varepsilon_{k-m+i,m}^{(n)})},\qquad k,n=0,1,\ldots,
\end{equation}
with the initial values
\begin{equation}
\varepsilon_{-m,m}^{(n)}=0, \ \varepsilon_{-m+1,m}^{(n)}=\varepsilon_{-m+2,m}^{(n)}=\cdots=\varepsilon_{-1,m}^{(n)}=n,\
\varepsilon_{0,m}^{(n)}=S_n, \qquad n=0, 1, \ldots \label{ini}
\end{equation}

Displaying these quantities in a double array similar to the $\varepsilon$--array, we see that this rule relates
$2m+2$ quantities located in an extended lozenge covering $m+2$ columns and two descending diagonals as showed below

$$
\begin{array}{ccccccc}
&\varepsilon_{k-m+1,m}^{(n)}&&&&&\\
\varepsilon_{k-m,m}^{(n+1)} && \varepsilon_{k-m+2,m}^{(n)}&&&&\\
&\varepsilon_{k-m+1,m}^{(n+1)}&&\ddots&&&\\
&&\ddots&&\ddots&&\\
&&&\ddots&&\varepsilon_{k,m}^{(n)}& \\
&&&&\varepsilon_{k-1,m}^{(n+1)}&&\varepsilon_{k+1,m}^{(n)}\\
&&&&&\varepsilon_{k,m}^{(n+1)}&
\end{array}
$$

The implementation of this algorithm using the technique of  ascending diagonals, as described above for the
$\varepsilon$--algorithm of Wynn, is more difficult, and it requires the storage of $m$ ascending
diagonals for computing the $(m+1)${\it th} one.

\vskip 2mm

In Section \ref{sse}, we will prove that, for all $k$ and $n$,  it holds
\begin{eqnarray}
 \varepsilon_{(m+1)k,m}^{(n)}&=&\frac{H_{k+1}(S_n)}{H_{k}(\Delta^{m+1} S_n)}, \label{slt1}\\
\varepsilon_{(m+1)(k-1)+1,m}^{(n)}&=&\frac{H_{k-1}(\Delta^{m+2}S_n)}{H_{k}(\Delta S_n)}, \label{slt2}\\
\varepsilon_{(m+1)(k-1)+i,m}^{(n)}&=&\frac{\Phi_{k+1}(\Delta^{i-1}S_n)}{H_{k}(\Delta^iS_n)},
\quad i=2, 3, \ldots, m,\label{slt3}
\end{eqnarray}
where the determinants $H_k$ and $\Phi_k$, which depend on $m$, are defined by
\begin{eqnarray*}
H_{k}(u_n)&=&\begin{vmatrix}
u_n             &  u_{n+1}           & \cdots  & u_{n+k-1}\\
\Delta^{m}u_n     & \Delta^{m}u_{n+1}    & \cdots  & \Delta^{m}u_{n+k-1}\\
\Delta^{2m}u_n     & \Delta^{2m}u_{n+1}    & \cdots  & \Delta^{2m}u_{n+k-1}\\
\vdots & \vdots & &\vdots\\
\Delta^{(k-1)m}u_n  & \Delta^{(k-1)m}u_{n+1} & \cdots  &
\Delta^{(k-1)m}u_{n+k-1}
\end{vmatrix},\quad
k=1, 2, \ldots ,\quad n=0, 1, \ldots,
\end{eqnarray*}
with $H_{-1}(u_n)=0$ and $\ H_{0}(u_n)=1$, and where

\begin{eqnarray*}
\Phi_{k}(u_n)&=&\begin{vmatrix}
n             &  n+1           & \cdots  & n+k-1\\
u_n             &  u_{n+1}           & \cdots  & u_{n+k-1}\\
\Delta^{m}u_n     & \Delta^{m}u_{n+1}    & \cdots  & \Delta^{m}u_{n+k-1}\\
\Delta^{2m}u_n     & \Delta^{2m}u_{n+1}    & \cdots  & \Delta^{2m}u_{n+k-1}\\
\vdots & \vdots & &\vdots\\
\Delta^{(k-2)m}u_n  & \Delta^{(k-2)m}u_{n+1} & \cdots  &
\Delta^{(k-2)m}u_{n+k-1}
\end{vmatrix},\quad
k=1, 2, \ldots, \quad n=0, 1, \ldots,
\end{eqnarray*}
with $\Phi_{-1}(u_n)=0$ and $\Phi_{0}(u_n)=1$.

Let us notice that, when $m=1$, $H_{k}(u_n)$ is identical to the usual Hankel determinant ${\cal H}_{k}(u_n)$.

\vskip 2mm

For proving these determinantal identities,
we will follow a procedure similar, although more difficult, to the procedure used in \cite{huweni}
(which is based on the Hirota's bilinear method) for
deriving a determinantal expression for a new acceleration algorithm obtained from the lattice Boussinesq equation.
However, instead of the Jacobi's determinantal identity, we will only use the Sylvester's one (which is,
in fact, the same after a permutation of rows and columns), and we will not use the Schwein's identity.

\vskip 2mm

Let us now define the {\it multistep Shanks' transformation} $e_{k,m} : (S_n) \longmapsto \{(e_{k,m}(S_n))\}$ by
\begin{equation}
e_{k,m}(S_n)=\varepsilon_{(m+1)k,m}^{(n)}=\frac{H_{k+1}(S_n)}{H_{k}(\Delta^{m+1} S_n)}, \qquad k,n=0,1,\ldots
\label{trf}
\end{equation}
Obviously, after proving (\ref{slt2}), we will also have
$$\varepsilon_{(m+1)k+1,m}^{(n)}=\frac{1}{e_{k,m}(\Delta S_n)}, \qquad k,n=0,1,\ldots,$$
a result similar to the second relation (\ref{se}) for the $\varepsilon$--algorithm of Wynn. Thus, only the quantities
$\varepsilon_{k,m}^{(n)}$'s whose first lower index is a multiple of $m+1$ are interesting for the purpose
of convergence acceleration. All the other ones are intermediate computations. The computation of
$e_{k,m}(S_n)=\varepsilon_{(m+1)k,m}^{(n)}$ needs the knowledge of $S_n,\ldots,S_{n+mk}$.

For simplicity, we will omit to indicate that all the symbols used in this paper depend on the fixed integer $m$.

\vskip 2mm

We see that, when $m=1$, the algorithm \eqref{gd} reduces to the $\varepsilon$--algorithm \eqref{epsi}, and the
transformation \eqref{trf} reduces to the Shanks' transformation \eqref{shanks}.
When $m=2$, the recursive rule \eqref{gd} reduces to the algorithm obtained
in \cite{huweni} from the lattice Boussinesq equation; see also \cite{Nij1,Nij2}.

\vskip 2mm

Let us mention that, due to \eqref{slt1}, the multistep Shanks' transformation can likewise be implemented by the $E$--algorithm \cite{ealg}
with $g_i(n)=\Delta^{im} S_n$ for $i=1,2,\ldots$, and for all $n$, and that we get, for all $k$ and $n$,
$E_k^{(n)}=e_{k,m}(S_n)$. Thus, by the fundamental property of the $E$--algorithm, the {\it kernel} of the
transformation \eqref{trf} (that is the set of sequences which are transformed into a constant sequence) is given by the
\begin{theorem}
~~\\
A necessary and sufficient condition that, for all $n$,
$e_{k,m}(S_n)=S$ is that there exist constants $a_1,\ldots,a_k$, $a_k \neq 0$, such that, for all $n$,
\[
S_n=S+a_1\Delta^{m}S_n+a_2\Delta^{2m}S_n+\cdots+a_{k}\Delta^{km}S_n.
\]
\end{theorem}

Let us remind that the kernel of the Shanks' transformation
$e_{km}: (S_n) \longmapsto (e_{km}(S_n)=\varepsilon_{2km}^{(n)})$
is the set of sequences such that, for all $n$,
$S_n=S+b_1 \Delta S_n+\cdots+b_{km} \Delta^{km} S_n$, where $b_1,\ldots,b_{km}$, $b_{km} \neq 0$, are
constants. Thus, we have the

\begin{corol}
~~\\
The kernel to the multistep Shanks' transformation $e_{k,m}$ is contained into the kernel of the Shanks'
transformation $e_{km}$.
\end{corol}

Moreover, due to the connection with the $E$--algorithm, all the convergence and acceleration results proved for it
\cite{ealg,ammp} also hold for the multistep Shanks' transformation.

\vskip 2mm

In the next Sections, we will link the multistep Shanks' transformation \eqref{trf} and the
multistep $\varepsilon$--algorithm \eqref{gd} by means of the Hirota's bilinear method.
First, in Section \ref{sshan}, some relations between the determinants ${H_{k}(\Delta^i S_n)}$ and
${\Phi_{k}(\Delta^i S_n)}$ will be established. We will only employ the Sylsvester's determinantal identity,
contrarily to the proofs given in \cite{wynn} and \cite{huweni} where the Schweins' determinantal identity is also used.
Then, Hirota's bilinear method will be presented in Section \ref{sshiro}. In Section \ref{sse}, we will show that the
quantities computed by the multistep $\varepsilon$--algorithm (\ref{gd}) are those defined in
the multistep Shanks' transformation (\ref{slt1})--(\ref{slt3}).
Conversely, in Section \ref{sss}, we will prove that the multistep Shanks' transformation (\ref{slt1})--(\ref{slt3})
can be implemented by the recursive rule (\ref{gd}) of the multistep $\varepsilon$--algorithm.

\section{Relations between determinants}
\label{sshan}

Let $A$ be a square matrix, $\alpha, \beta, \gamma$
and $\delta$ numbers, $a,b,c$ and $d$ vectors of the same dimension as $A$. Let $M$ be the matrix
$$M=\left(
\begin{array}{ccc}
\alpha & a^T & \beta \\
b & A & c \\
\gamma & d^T & \delta
\end{array}
\right).$$
The {\it Sylvester's determinantal identity} is
$$|M|\cdot|A|=\left|
\begin{array}{cc}
\alpha & a^T \\
b & A \\
\end{array}
\right|\cdot
\left|
\begin{array}{cc}
A & c \\
d^T & \delta
\end{array}
\right|-
\left|
\begin{array}{cc}
a^T & \beta \\
A & c \\
\end{array}
\right|\cdot
\left|
\begin{array}{cc}
b & A \\
\gamma & d^T
\end{array}
\right|.$$

Let us now prove some determinantal identities that will be useful in the sequel.

\vskip 2mm
\begin{lemma}
\begin{eqnarray}
H_{k+1}(\Delta S_n)H_{k}(\Delta^{m}S_{n+1})\!\!\!\!&=&\!\!\!\!H_{k}(\Delta^{m+1}S_n)H_{k+1}(S_{n+1})-
H_{k}(\Delta^{m+1}S_{n+1})H_{k+1}(S_{n}).  \label{id3}
\end{eqnarray}
\end{lemma}

\noindent {\bf Proof:}
we consider the determinant
\begin{equation*}
D_1=\left|
\begin{array}{cccc}
1  &   1  &\cdots &   1  \\
\Delta^{i+m}S_n & \Delta^{i+m}S_{n+1} & \cdots & \Delta^{i+m}S_{n+k+1} \\
\vdots & \vdots && \vdots \\
\Delta^{i+km}S_n & \Delta^{i+km}S_{n+1} & \cdots & \Delta^{i+km}S_{n+k+1}\\
\Delta^{i}S_n & \Delta^{i}S_{n+1} & \cdots & \Delta^{i}S_{n+k+1}
\end{array}
\right|=(-1)^kH_{k+1}(\Delta^{i+1}S_n).
\end{equation*}
The second expression for $D_1$ is obtained by replacing each column, from the last one, by its difference with the
previous one. Thus, we get a determinant whose first row
only contains 0 except in the first column where the element is equal to 1. Expanding this determinant with respect
to its first row, and putting its last row as the first one, we see that $D_1=(-1)^kH_{k+1}(\Delta^{i+1}S_n)$. Let us now
apply the Sylvester's identity to the first expression of $D_1$, and perform a similar manipulation on the rows
and the columns of the other determinants, we obtain
\begin{equation}
H_{k+1}(\Delta^{i+1}S_n)H_{k}(\Delta^{i+m}S_{n+1})\!\!=\!\!H_{k}(\Delta^{i+m+1}S_n)H_{k+1}(\Delta^{i}S_{n+1})-
H_{k}(\Delta^{i+m+1}S_{n+1})H_{k+1}(\Delta^{i}S_{n}).\nonumber
\end{equation}
Setting $i=0$ in this relation, we get \eqref{id3}. \fin

\vskip 1mm

A similar identity, which will be used in the sequel, also holds if $S_n$ is replaced by $\Delta S_n$.

\vskip 2mm
\begin{lemma}
\begin{equation}
H_{k}(\Delta^{i+1}S_n)H_{k-1}(\Delta^{i}S_{n+1})=H_{k-1}(\Delta^{i+1}S_n)H_{k}(\Delta^{i}S_{n+1})-
H_{k-1}(\Delta^{i+1}S_{n+1})H_{k}(\Delta^{i}S_{n}). \label{id2}
\end{equation}
\end{lemma}

\noindent {\bf Proof:}
let $D_2$ be the determinant obtained from $D_1$ by replacing $k$ by $k-1$, and moving the last row
to the second position.
Replacing each column, from the last one, by its difference with the previous one, we see
that $D_2=H_{k}(\Delta^{i+1}S_n)$, and, applying the Sylvester's identity to it, we get \eqref{id2}. \fin

\vskip 2mm
\begin{lemma}
\begin{eqnarray}
H_{k}(\Delta S_{n})H_{k}(\Delta^mS_{n+1})\!\!\!\!&=&\!\!\!\!H_{k}(\Delta^{m+1}S_n)H_{k}(S_{n+1})-
H_{k+1}(S_n)H_{k-1}(\Delta^{m+1}S_{n+1}).  \label{id7}
\end{eqnarray}
\end{lemma}

\noindent {\bf Proof:}
Setting $i=m$ in \eqref{id2}, we have
\begin{equation}
H_{k}(\Delta^{m+1}S_n)H_{k-1}(\Delta^{m}S_{n+1})=H_{k-1}(\Delta^{m+1}S_n)H_{k}(\Delta^{m}S_{n+1})-
H_{k-1}(\Delta^{m+1}S_{n+1})H_{k}(\Delta^{m}S_{n}).\label{eq8.1}
\end{equation}

Applying now the Sylvester's identity to the determinant $H_{k+1}(\Delta^iS_n)$, we get
\begin{equation}
H_{k+1}(\Delta^iS_n)H_{k-1}(\Delta^{i+m}S_{n+1})=H_{k}(\Delta^iS_n)H_{k}(\Delta^{i+m}S_{n+1})-
H_{k}(\Delta^iS_{n+1})H_{k}(\Delta^{i+m}S_{n}). \nonumber
\end{equation}

Setting $i=0$ in this relation, we obtain
\begin{equation}
H_{k+1}(S_n)H_{k-1}(\Delta^{m}S_{n+1})=H_{k}(S_n)H_{k}(\Delta^{m}S_{n+1})-H_{k}(S_{n+1})H_{k}(\Delta^{m}S_{n}).
\label{eq8.2}
\end{equation}

Then, we multiply (\ref{eq8.1}) by $H_{k}(S_{n+1})$, we multiply (\ref{eq8.2}) by $H_{k-1}(\Delta^{m+1}S_{n+1})$,
and we subtract. It gives
\begin{eqnarray*}
H_{k-1}(\Delta^{m}S_{n+1})[H_{k}(\Delta^{m+1}S_{n})H_{k}(S_{n+1})-H_{k+1}(S_n)H_{k-1}
(\Delta^{m+1}S_{n+1})]&&\\
=H_{k}(\Delta^{m}S_{n+1})[H_{k-1}(\Delta^{m+1}S_{n})H_{k}(S_{n+1})-H_{k}(S_{n})H_{k-1}
(\Delta^{m+1}S_{n+1})].&&
\end{eqnarray*}
Using (\ref{id3}), we see that the bracket in the right hand side is equal to $H_{k}(\Delta S_n)H_{k-1}(\Delta^mS_{n+1})$.
After simplifying both sides by $H_{k-1}(\Delta^mS_{n+1})$, we obtain \eqref{id7}. \fin

\vskip 1mm

A similar identity, which will be used in the sequel, also holds if $S_n$ is replaced by $\Delta S_n$.

\vskip 2mm

\begin{lemma}
\begin{equation}
H_{k}(\Delta^iS_{n+1})H_{k-1}(\Delta^{i+2}S_{n})=H_{k}(\Delta^{i+1}S_{n})\Phi_{k}(\Delta^iS_{n+1})-
H_{k-1}(\Delta^{i+1}S_{n+1})\Phi_{k+1}(\Delta^iS_n).\label{id9}
\end{equation}
\end{lemma}

\noindent {\bf Proof:} we consider the determinant
\begin{equation*}
D_3=\left|
\begin{array}{cccc}
1 & 1 & \cdots & 1 \\
\Delta^{i}S_n & \Delta^{i}S_{n+1} & \cdots & \Delta^{i}S_{n+k+1} \\
\Delta^{i+m}S_n & \Delta^{i+m}S_{n+1} & \cdots & \Delta^{i+m}S_{n+k+1} \\
\vdots & \vdots && \vdots \\
\Delta^{i+(k-1)m}S_n & \Delta^{i+(k-1)m}S_{n+1} & \cdots & \Delta^{i+(k-1)m}S_{n+k+1}\\
n & n+1 & \cdots & n+k+1
\end{array}
\right|.
\end{equation*}
Obviously, we also have
\begin{eqnarray*}
D_3&=&(-1)^k\left|
\begin{array}{cccc}
1 & 1 & \cdots & 1 \\
n & n+1 & \cdots & n+k+1 \\
\Delta^{i}S_n & \Delta^{i}S_{n+1} & \cdots & \Delta^{i}S_{n+k+1} \\
\vdots & \vdots && \vdots \\
\Delta^{i+(k-1)m}S_n & \Delta^{i+(k-1)m}S_{n+1} & \cdots & \Delta^{i+(k-1)m}S_{n+k+1}
\end{array}
\right|\\
&=&(-1)^k\left|
\begin{array}{cccc}
1 & 1 & \cdots & 1 \\
\Delta^{i+1}S_n & \Delta^{i+1}S_{n+1} & \cdots & \Delta^{i+1}S_{n+k} \\
\vdots & \vdots && \vdots \\
\Delta^{i+1+(k-1)m}S_n & \Delta^{i+1+(k-1)m}S_{n+1} & \cdots & \Delta^{i+1+(k-1)m}S_{n+k}
\end{array}
\right|\\
&=&(-1)^kH_{k}(\Delta^{i+2}S_n).
\end{eqnarray*}
We apply now the Sylvester's identity to the first expression of $D_3$ given above, and replace $D_3$ by
$(-1)^kH_{k}(\Delta^{i+2}S_n)$. We get, after similar manipulations on the columns of the other determinants,
\begin{equation}
H_{k}(\Delta^{i+2}S_n)H_{k}(\Delta^{i}S_{n+1})=H_{k}(\Delta^{i+1}S_n)\Phi_{k+1}(\Delta^{i}S_{n+1})-
H_{k}(\Delta^{i+1}S_{n+1})\Phi_{k+1}(\Delta^{i}S_n). \label{id6}
\end{equation}

Then, we apply the Sylvester's identity to the determinant $\Phi_{k+1}(\Delta^i S_n)$. We get
\begin{equation}
\Phi_{k+1}(\Delta^iS_n)H_{k-1}(\Delta^{i}S_{n+1})=\Phi_{k}(\Delta^i S_n)H_{k}(\Delta^{i}S_{n+1})-
\Phi_{k}(\Delta^{i}S_{n+1})H_{k}(\Delta^{i}S_{n}). \nonumber
\end{equation}

We multiply this identity by $H_{k-1}(\Delta^{i+1}S_{n+1})$, we multiply (\ref{id2}) by $\Phi_{k}(\Delta^{i}S_{n+1})$,
and we subtract. It gives
\begin{eqnarray*}
H_{k-1}(\Delta^{i}S_{n+1})[\Phi_{k+1}(\Delta^{i}S_n)H_{k-1}(\Delta^{i+1}S_{n+1})-H_{k}(\Delta^{i+1}S_n)
\Phi_{k}(\Delta^{i}S_{n+1})]&&\\
=H_{k}(\Delta^{i}S_{n+1})[\Phi_{k}(\Delta^{i}S_n)H_{k-1}(\Delta^{i+1}S_{n+1})-H_{k-1}(\Delta^{i+1}S_{n})
\Phi_{k}(\Delta^{i}S_{n+1})].&&
\end{eqnarray*}
Using (\ref{id6}), we see that the bracket in the right hand side is equal to $-H_{k-1}(\Delta^{i+2}S_n)
H_{k-1}(\Delta^iS_{n+1})$.
After simplifying both sides by $H_{k-1}(\Delta^iS_{n+1})$, we obtain \eqref{id9}. \fin

\vskip 2mm
\begin{lemma}
\begin{equation}
H_{k}(\Delta S_n)H_{k-2}(\Delta^{m+1}S_{n+1})=H_{k-1}(\Delta^{m+1}S_{n+1})H_{k-1}(\Delta S_{n})
-H_{k-1}(\Delta^{m+1}S_{n})H_{k-1}(\Delta S_{n+1}). \label{cs1}
\end{equation}
\end{lemma}

\noindent {\bf Proof:} we consider the determinant
\begin{equation*}
D_4=\left|
\begin{array}{cccc}
S_n & S_{n+1} & \cdots & S_{n+k} \\
1  &   1  &\cdots &   1  \\
\Delta^{m}S_n & \Delta^{m}S_{n+1} & \cdots & \Delta^{m}S_{n+k} \\
\vdots & \vdots && \vdots \\
\Delta^{(k-1)m}S_n & \Delta^{(k-1)m}S_{n+1} & \cdots & \Delta^{(k-1)m}S_{n+k}\\
\end{array}
\right|.
\end{equation*}
After exchanging the first row and the second row, we see that $D_4=-H_{k}(\Delta S_n)$. Let us now
apply the Sylvester's identity to $D_4$, and perform a similar manipulation on the first and second row of the various
determinants. We obtain \eqref{cs1}. \fin

\section{The Hirota's bilinear method}
\label{sshiro}

The Hirota's bilinear method \cite{hiro} is a technique which could be much useful for solving certain nonlinear differential
and difference equations. It consists in expressing the unknown as a ratio and, then, in treating separately
the numerator and the denominator.

We will now apply this method to the multistep $\varepsilon$--algorithm, and set
\begin{equation}
\varepsilon_{k,m}^{(n)}=\frac{G_{k}^n}{F_{k}^n}. \label{hiro}
\end{equation}

 We first have the
 \begin{lemma}
 \label{lem1}
 \begin{equation}
(F_{k+m+1}^{n}G_{k}^{n+1}-F_{k}^{n+1}G_{k+m+1}^{n})
\prod_{i=1}^m(F_{k+i}^{n}G_{k+i}^{n+1}-F_{k+i}^{n+1}G_{k+i}^{n})
=-F_{k+1}^{n}F_{k+m}^{n+1}
\prod_{i=1}^mF_{k+i+1}^{n}F_{k+i-1}^{n+1}. \label{rre1}
\end{equation}
 \end{lemma}

\noindent {\bf Proof:}
~~\\
Plugging (\ref{hiro}) into the recursive rule (\ref{gd}) of the $\varepsilon$--algorithm, we get
\begin{eqnarray}
\frac{G_{k+1}^n}{F_{k+1}^n}-\frac{G_{k-m}^{n+1}}{F_{k-m}^{n+1}}&=&\displaystyle \frac{1}
{\prod_{i=1}^m\left(\displaystyle \frac{G_{k-m+i}^{n+1}}{F_{k-m+i}^{n+1}}-\frac{G_{k-m+i}^{n}}{F_{k-m+i}^{n}}\right)}
\nonumber \\
\frac{F_{k-m}^{n+1}G_{k+1}^{n}-F_{k+1}^{n}G_{k-m}^{n+1}}{F_{k+1}^{n}F_{k-m}^{n+1}}&=&
\frac{\prod_{i=1}^mF_{k-m+i}^{n}F_{k-m+i}^{n+1}}{\prod_{i=1}^m(F_{k-m+i}^{n}G_{k-m+i}^{n+1}
-F_{k-m+i}^{n+1}G_{k-m+i}^{n})}. \label{rr}
\end{eqnarray}

Now, we cross--multiply the numerator of one side by the denominator of the other side, and we equate both sides.
Replacing $k$ by $k+m$ and changing the sign, the equation \eqref{rr} becomes \eqref{rre1}
since
$$F_{k+m+1}^{n}F_{k}^{n+1}
\prod_{i=1}^mF_{k+i}^{n}F_{k+i}^{n+1}
=\prod_{i=1}^{m+1}F_{k+i}^{n}\prod_{i=0}^mF_{k+i}^{n+1}. \fin$$

The second preliminary result is contained in the

\begin{lemma}
~~\\
\label{lem2}
If the following relations hold
\begin{eqnarray}
&&F_{(m+1)k+1}^{n}G_{(m+1)k+1}^{n+1}-F_{(m+1)k+1}^{n+1}G_{(m+1)k+1}^{n}=-F_{(m+1)k+2}^{n}
F_{(m+1)k}^{n+1},\label{bt2}\\
&&F_{(m+1)k+1}^{n}G_{(m+1)(k-1)+1}^{n+1}-F_{(m+1)(k-1)+1}^{n+1}G_{(m+1)k+1}^{n}
=-F_{(m+1)(k-1)+2}^{n}F_{(m+1)k}^{n+1},\label{bt4}
\end{eqnarray}
and,  for $i=2,\ldots,m+1$,
\begin{eqnarray}
&&F_{(m+1)k+i}^{n}G_{(m+1)k+i}^{n+1}-F_{(m+1)k+i}^{n+1}G_{(m+1)k+i}^{n}=F_{(m+1)k+i+1}^{n}
F_{(m+1)k+i-1}^{n+1},\label{bt1}\\
&&F_{(m+1)k+i}^{n}G_{(m+1)(k-1)+i}^{n+1}-F_{(m+1)(k-1)+i}^{n+1}G_{(m+1)k+i}^{n}
=F_{(m+1)(k-1)+i+1}^{n}F_{(m+1)k+i-1}^{n+1},\label{bt3}
\end{eqnarray}
then \eqref{rre1} follows.
\end{lemma}

\noindent {\bf Proof:}
~~\\
Let us first notice that, taking $i=1$ in \eqref{bt1} and \eqref{bt3} gives \eqref{bt2} and \eqref{bt4}, respectively,
after a change in the signs of their right hand sides.

Let us separate \eqref{rre1} into the product of two relations, and prove that each of the following
formulae holds separately
\begin{equation}
\prod_{i=1}^m(F_{k+i}^{n}G_{k+i}^{n+1}-F_{k+i}^{n+1}G_{k+i}^{n})=\pm \prod_{i=1}^mF_{k+i+1}^{n}F_{k+i-1}^{n+1}, \label{p1}
\end{equation}
which are the products appearing in both sides of \eqref{rre1}, and
\begin{equation}
F_{k+m+1}^{n}G_{k}^{n+1}-F_{k}^{n+1}G_{k+m+1}^{n}=\mp F_{k+1}^{n}F_{k+m}^{n+1}, \label{p2}
\end{equation}
which are its remaining parts. Then, multiplying together \eqref{p1} and \eqref{p2}, we will obtain
\eqref{rre1}, but we must notice that the signs used in \eqref{p1} and \eqref{p2} have to be opposite.

\vskip 1mm

Let us assume that (\ref{bt2})--(\ref{bt3}) hold true.
The proofs of the relations \eqref{p1} and \eqref{p2} have to be separated into three cases according to the value of
$k$ in \eqref{rre1}.

\vskip 1mm

\noindent $\bullet$ $k$ replaced by $(m+1)k$ in \eqref{p1} and \eqref{p2}.

$\circ$ Multiplying together the relations \eqref{bt1} for $i=2,\ldots,m$, and then multiplying each of its sides by the
corresponding side of \eqref{bt2} (which brings a change in the sign) proves \eqref{p1}, with the sign $-$, when
$k$ is replaced by $(m+1)k$ in \eqref{p1}.

$\circ$ Replacing $k$ by $(m+1)k$ in \eqref{p2}, with the sign $+$, gives
$$F_{(m+1)(k+1)}^{n}G_{(m+1)k}^{n+1}-F_{(m+1)k}^{n+1}G_{(m+1)(k+1)}^{n}=F_{(m+1)k+1}^{n}F_{(m+1)k+m}^{n+1},$$
which is \eqref{bt3} for $i=m+1$.

$\circ$ We get \eqref{rre1} by multiplying together the two relations.

\vskip 1mm

\noindent $\bullet$ $k$ replaced by $(m+1)k+1$ in \eqref{p1} and \eqref{p2}.

$\circ$ Multiplying together the relations \eqref{bt1} for $i=2,\ldots,m+1$. The result is the same as adding 1 to all the
lower indexes, and making the product for $i=1,\ldots,m$, which is \eqref{p1} with the sign $+$.

$\circ$ Replace $k$ by $(m+1)k+1$ in \eqref{p2}, with the sign $-$. It is exactly \eqref{bt4} with $k+1$ instead of $k$.

$\circ$ Multiplying together the two relations, we obtain \eqref{rre1}.

\vskip 1mm

\noindent $\bullet$  $k$ replaced by $(m+1)k+j$, for $j=2,\ldots,m$, in \eqref{p1} and \eqref{p2}.

$\circ$ Let $2 \leq j \leq m$ be fixed. Let us write that \eqref{bt1} holds with $2 \leq i+j \leq m+1$ instead of $i$,
that is for $i=1,\ldots,m-j+1$,
$$F_{(m+1)k+i+j}^{n}G_{(m+1)k+i+j}^{n+1}-F_{(m+1)k+i+j}^{n+1}G_{(m+1)k+i+j}^{n}=F_{(m+1)k+i+j+1}^{n}
F_{(m+1)k+i+j-1}^{n+1}.$$
Multiply together these relations for $i=1,\ldots,m-j+1$, and, then, make their product
for $i=m-j+2,\ldots,m$. When $i=m-j+2$, the left hand side of this expression becomes
$F_{(m+1)k+m+2}^{n}G_{(m+1)k+m+2}^{n+1}-F_{(m+1)k+m+2}^{n+1}G_{(m+1)k+m+2}^{n}$, and its right hand side is
equal to $F_{(m+1)k+m+3}^{n}F_{(m+1)k+m+1}^{n+1}$, that is, respectively,
$F_{(m+1)(k+1)+1}^{n}G_{(m+1)(k+1)+1}^{n+1}-F_{(m+1)(k+1)+1}^{n+1}G_{(m+1)(k+1)+1}^{n}$, and
$F_{(m+1)(k+1)+2}^{n} F_{(m+1)(k+1)}^{n+1}$.
Thus, by \eqref{bt2} with $k$ replaced by $k+1$, these two expressions are equal after changing the sign in one side.
For $i=m-j+3$, we have
$$F_{(m+1)k+m+3}^{n}G_{(m+1)k+m+3}^{n+1}-F_{(m+1)k+m+3}^{n+1}G_{(m+1)k+m+3}^{n}=F_{(m+1)k+m+4}^{n}
F_{(m+1)k+m+2}^{n+1},$$
that is
$$F_{(m+1)(k+1)+2}^{n}G_{(m+1)(k+1)+2}^{n+1}-F_{(m+1)(k+1)+2}^{n+1}G_{(m+1)(k+1)+2}^{n}=F_{(m+1)(k+1)+3}^{n}
F_{(m+1)(k+1)+1}^{n+1},$$
which is \eqref{bt1} with $k+1$ instead of $k$. And so on until $i=m$.
Thus in the products from $i=m-j+2$ to $m$, the sign is changed in one, and only one, of the expressions
due to \eqref{bt2}, and we finally obtain \eqref{p1} with the sign $-$.

$\circ$ Let us replace $k$ by $(m+1)k+j$ in \eqref{p2}, with the sign $+$. We get
$$F_{(m+1)(k+1)+j}^{n}G_{(m+1)k+j}^{n+1}-F_{(m+1)k+j}^{n+1}G_{(m+1)(k+1)+j}^{n}=
F_{(m+1)k+j}^{n}F_{(m+1)(k+1)+j-1}^{n+1},$$
which is \eqref{bt3} when $k$ is replaced by $k+1$.

$\circ$ The product of the two relations gives \eqref{rre1}.

\vskip 1mm

Thus, \eqref{rre1} have now been proved for all values of $k$. \fin

\vskip 2mm

Finally, we are able to prove the

\begin{theorem}
~~\\
\label{lem3}
The relation \eqref{rre1} holds with the $F_k^n$'s and the $G_k^n$'s given by the following relations,
for $k=0,1,\ldots$,
\begin{eqnarray}
F_{(m+1)(k-1)+i}^n&=&H_{k}(\Delta^iS_n),  \quad i=1, 2, \ldots, m+1\label{bf1}\\
G_{(m+1)(k-1)+1}^n&=&H_{k-1}(\Delta^{m+2}S_n), \quad G_{(m+1)k}^n=H_{k+1}(S_n), \label{bf11} \\
G_{(m+1)(k-1)+i}^n&=&\Phi_{k+1}(\Delta^{i-1}S_n), \quad i=2, 3, \ldots, m.\label{bf2}
\end{eqnarray}
\end{theorem}
\noindent {\bf Proof:}
~~\\
We are now able to prove \eqref{bt2}--\eqref{bt3}, with the $F_k^n$'s and $G_k^n$'s given by \eqref{bf1}--\eqref{bf2}.
Replacing the determinants in (\ref{id3}) by their expressions, we obtain
\begin{equation*}
F_{(m+1)(k+1)}^{n}G_{(m+1)(k+1)}^{n+1}-F_{(m+1)(k+1)}^{n+1}G_{(m+1)(k+1)}^{n}=
F_{(m+1)(k+1)+1}^{n}F_{(m+1)k+m}^{n+1}
\end{equation*}
which corresponds to \eqref{bt1} for the case $i=m+1$. Replacing the determinants in (\ref{id6})
by their expressions, we obtain the bilinear equation \eqref{bt1} for the cases $i=2, 3, \ldots, m$,
which completes the proof of the equation \eqref{bt1}.

Replacing the determinants in (\ref{id3}) and (\ref{id7}), both with $\Delta S_n$
instead of $S_n$, by their expressions, we see that the equations \eqref{bt2}
and \eqref{bt4} are satisfied.

Then, replacing the determinants in (\ref{id7}) by their expressions, we obtain
\begin{equation*}
F_{(m+1)k}^{n}G_{(m+1)(k-1)}^{n+1}-F_{(m+1)(k-1)}^{n+1}G_{(m+1)k}^{n}=F_{(m+1)(k-1)+1}^{n}F_{(m+1)(k-1)+m}^{n+1}
\end{equation*}
which corresponds to \eqref{bt3} for the case $i=m+1$, while replacing the determinants in (\ref{id9})
by their expressions, we get the bilinear equation \eqref{bt3} for $i=2, 3, \ldots, m$,
which completes the proof for the equation \eqref{bt3}.

Since the identities \eqref{bt2}--\eqref{bt3} hold, then \eqref{rre1} follows with the $F_k^n$'s and the $G_k^n$'s
given by \eqref{bf1}-- \eqref{bf2}. \fin

\vskip 2mm

We also have the

\begin{corol}
\begin{equation}
F_{k+m+1}^nF_{k-1}^{n+1}=F_{k}^nF_{k+m}^{n+1}-F_{k+m}^nF_{k}^{n+1}. \nonumber
\end{equation}
\end{corol}
\noindent {\bf Proof:}
~~\\
Replacing the determinants in \eqref{id2} by their expressions given by \eqref{bf1}, and $k$ by $k+1$,
we obtain, for $i=1,\ldots,m$, the following relation without any $G_k^n$

\begin{equation}
F_{(m+1)k+i+1}^nF_{(m+1)(k-1)+i}^{n+1}=F_{(m+1)(k-1)+i+1}^nF_{(m+1)k+i}^{n+1}
-F_{(m+1)k+i}^nF_{(m+1)(k-1)+i+1}^{n+1}. \nonumber
\end{equation}

Similarly, the determinantal identity \eqref{cs1} leads, after replacing $k$ by $k+2$, to
\begin{eqnarray}
F_{(m+1)(k+1)+1}^nF_{(m+1)k}^{n+1}=F_{(m+1)k+1}^nF_{(m+1)(k+1)}^{n+1}
-F_{(m+1)(k+1)}^nF_{(m+1)k+1}^{n+1}, \nonumber
\end{eqnarray}
which is the preceding relation for $i=m+1$. Thus, changing $mk$ into $m$, these two identities can be gathered into the single formula of the Corollary. \fin

\section{From the multistep $\varepsilon$--algorithm to the multistep Shanks' transformation}
\label{sse}

By comparing \eqref{slt1}--\eqref{slt3} with the determinantal formulae \eqref{bf1}--\eqref{bf2} of the Theorem \ref{lem3} for the $F_k^n$'s and the $G_k^n$'s issued from the Hirota's method,
we are now able to give the determinantal formulae for the multistep $\varepsilon$--algorithm.
Consequently, from the Lemmas \ref{lem1}, \ref{lem2}, and the Theorem \ref{lem3}, we have the

\begin{theorem}
~~\\
The quantities $\varepsilon_{k,m}^{(n)}$ computed by the multistep $\varepsilon$--algorithm \eqref{gd}, with the initializations \eqref{ini}, are expressed by the ratios of determinants \eqref{slt1}, \eqref{slt2}, and \eqref{slt3}.
\end{theorem}

\vskip 2mm

Thus, starting from the determinantal identities between ${H_{k}(\Delta^iS_n)}$ and ${\Phi_{k}(\Delta^iS_n)}$,
we proved that \eqref{bt2}--\eqref{bt3} are satisfied with the determinantal formulae \eqref{bf1}--\eqref{bf2} for
the $F_k^n$'s and the $G_k^n$'s. Then, \eqref{rre1} followed, and we concluded that the determinantal
expressions \eqref{slt1}--\eqref{slt3} for the $\varepsilon_{k,m}^{(n)}$'s hold true.
Notice that all these results were obtained without using the rule (\ref{gd}) of the multistep $\varepsilon$--algorithm.

\vskip 2mm

Let us remind that, as noticed in \cite{5benc} and fully explained in \cite{5fdf}, we have
$$\varepsilon_{(m+1)k,m}^{(n)}=\frac{f_{k,m}(S_n,\ldots,S_{n+(m+1)k})}{Df_{k,m}(S_n,\ldots,S_{n+(m+1)k})}, \qquad
\varepsilon_{(m+1)k+1,m}^{(n)}=\frac{Df_{k,m}(\Delta S_n,\ldots,\Delta S_{n+(m+1)k})}{f_{k,m}(\Delta S_n,\ldots,
\Delta S_{n+(m+1)k})},$$
where $f_{k,m}$ is a function depending on $(m+1)k+1$ variables and such that $D^2f_{k,m} \equiv 0$, where $Df_{k,m}$
denotes
the sum of the partial derivatives of $f_{k,m}$. Thus, we obtain the following connection with Hirota's bilinear method
$$
\begin{array}{ll}
G_{(m+1)k}^n=f_{k,m}(S_n,\ldots,S_{n+(m+1)k}), &F_{(m+1)k}^n=Df_{k,m}(S_n,\ldots,S_{n+(m+1)k}), \\
G_{(m+1)k+1}^n=Df_{k,m}(\Delta S_n,\ldots,\Delta S_{n+(m+1)k}), &F_{(m+1)k+1}^n=f_{k,m}(\Delta S_n,\ldots,
\Delta S_{n+(m+1)k}),
\end{array}
$$
and, according to this theory, the multistep Shanks' transformation is quasilinear that is
$e_{k,m}(aS_n+b)=ae_{k,m}(S_n)+b$, a result which can be seen directly from \eqref{trf}.

\section{From the multistep Shanks' transformation to the multistep $\varepsilon$--algorithm}
\label{sss}

We will show now how to derive the recursive rule (\ref{gd}) of the multistep $\varepsilon$--algorithm from the
definition (\ref{slt1})--\eqref{slt3} of the multistep Shanks' transformation.

 From the determinantal identity \eqref{id7}, we get
\begin{eqnarray}
\varepsilon_{(m+1)(k+1),m}^{(n)}-\varepsilon_{(m+1)k,m}^{(n+1)}&=&\frac{H_{k+2}(S_n)}{H_{k+1}(\Delta^{m+1}S_n)}-
\frac{H_{k+1}(S_{n+1})}{H_{k}(\Delta^{m+1}S_{n+1})}\nonumber\\
&=&\frac{H_{k+2}(S_n)H_{k}(\Delta^{m+1}S_{n+1})-H_{k+1}(S_{n+1})H_{k+1}(\Delta^{m+1}S_n)}{H_{k+1}
(\Delta^{m+1}S_n)H_{k}(\Delta^{m+1}S_{n+1})}\nonumber\\
&=&-\frac{H_{k+1}(\Delta S_n)H_{k+1}(\Delta^mS_{n+1})}{H_{k+1}(\Delta^{m+1}S_n)H_{k}
(\Delta^{m+1}S_{n+1})}.\label{l1}
\end{eqnarray}

Similarly, by the identity \eqref{id7} with $S_n$ replaced by $\Delta S_n$, we get
\begin{eqnarray}
\varepsilon_{(m+1)(k+1)+1,m}^{(n)}-\varepsilon_{(m+1)k+1,m}^{(n+1)}&=&\frac{H_{k+1}(\Delta^{m+2}S_n)}{H_{k+2}
(\Delta S_n)}-\frac{H_{k}(\Delta^{m+2}S_{n+1})}{H_{k+1}(\Delta S_{n+1})}\nonumber\\
&=&\frac{H_{k+1}(\Delta^{m+2}S_n)H_{k+1}(\Delta S_{n+1})-H_{k}(\Delta^{m+2}S_{n+1})H_{k+2}(\Delta S_n)}
{H_{k+2}(\Delta S_n)H_{k+1}(\Delta S_{n+1})}\nonumber\\
&=&\frac{H_{k+1}(\Delta^2 S_n)H_{k+1}(\Delta^{m+1}S_{n+1})}{H_{k+2}(\Delta S_n)H_{k+1}
(\Delta S_{n+1})}.\label{l2}
\end{eqnarray}

We also get the following relation from the identity \eqref{id9}
\begin{eqnarray}
\varepsilon_{(m+1)(k+1)+i,m}^{(n)}-\varepsilon_{(m+1)k+i,m}^{(n+1)}&=&\frac{\Phi_{k+3}(\Delta^{i-1}S_n)}
{H_{k+2}(\Delta^i S_n)}-\frac{\Phi_{k+2}(\Delta^{i-1}S_{n+1})}{H_{k+1}(\Delta^i S_{n+1})}\nonumber\\
&=&\frac{\Phi_{k+3}(\Delta^{i-1}S_n)H_{k+1}(\Delta^i S_{n+1})-\Phi_{k+2}(\Delta^{i-1}S_{n+1})H_{k+2}
(\Delta^i S_n)}{H_{k+2}(\Delta^i S_n)H_{k+1}(\Delta^i S_{n+1})}\nonumber\\
&=&-\frac{H_{k+2}(\Delta^{i-1}S_{n+1})H_{k+1}(\Delta^{i+1}S_{n})}{H_{k+2}(\Delta^i S_n)H_{k+1}
(\Delta^i S_{n+1})}, \quad i=2, 3, \ldots, m. \label{l3}
\end{eqnarray}

Besides, from the identity \eqref{id3}, with $\Delta S_n$ instead of $S_n$, we get
\begin{eqnarray}
\varepsilon_{(m+1)k+1,m}^{(n+1)}-\varepsilon_{(m+1)k+1,m}^{(n)}&=&\frac{H_{k}(\Delta^{m+2}S_{n+1})}{H_{k+1}
(\Delta S_{n+1})}-\frac{H_{k}(\Delta^{m+2}S_{n})}{H_{k+1}(\Delta S_{n})}\nonumber\\
&=&\frac{H_{k}(\Delta^{m+2}S_{n+1})H_{k+1}(\Delta S_{n})-H_{k}(\Delta^{m+2}S_{n})H_{k+1}
(\Delta S_{n+1})}{H_{k+1}(\Delta S_{n+1})H_{k+1}(\Delta S_{n})}\nonumber\\
&=&-\frac{H_{k+1}(\Delta^2 S_n)H_{k}(\Delta^{m+1}S_{n+1})}{H_{k+1}(\Delta S_{n+1})H_{k+1}
(\Delta S_{n})}. \label{h1}
\end{eqnarray}

 From \eqref{id3}, we get
\begin{eqnarray}
\varepsilon_{(m+1)(k+1),m}^{(n+1)}-\varepsilon_{(m+1)(k+1),m}^{(n)}&=&\frac{H_{k+2}(S_{n+1})}{H_{k+1}
(\Delta^{m+1} S_{n+1})}-\frac{H_{k+2}(S_{n})}{H_{k+1}(\Delta^{m+1}S_{n})}\nonumber\\
&=&\frac{H_{k+2}(S_{n+1})H_{k+1}(\Delta^{m+1} S_{n})-H_{k+2}(S_{n})H_{k+1}(\Delta^{m+1}S_n)}
{H_{k+1}(\Delta^{m+1}S_{n+1})H_{k+1}(\Delta^{m+1}S_{n})}\nonumber\\
&=&\frac{H_{k+2}(\Delta S_n)H_{k+1}(\Delta^{m}S_{n+1})}{H_{k+1}(\Delta^{m+1}S_{n+1})
H_{k+1}(\Delta^{m+1}S_{n})}.\label{h2}
\end{eqnarray}

 Finally, from the identity \eqref{id6}, we have, for $i=2,\ldots,m$,
\begin{eqnarray}
\varepsilon_{(m+1)k+i,m}^{(n+1)}-\varepsilon_{(m+1)k+i,m}^{(n)}&=&\frac{\Phi_{k+2}(\Delta^{i-1}S_{n+1})}
{H_{k+1}(\Delta^i S_{n+1})}-\frac{\Phi_{k+2}(\Delta^{i-1}S_n)}{H_{k+1}(\Delta^i S_{n})}\nonumber\\
&=&\frac{\Phi_{k+2}(\Delta^{i-1}S_{n+1})H_{k+1}(\Delta^i S_{n})-\Phi_{k+2}(\Delta^{i-1}S_n)
H_{k+1}(\Delta^i S_{n+1})}{H_{k+1}(\Delta^i S_{n+1})H_{k+1}(\Delta^iS_{n})}\nonumber\\
&=&\frac{H_{k+1}(\Delta^{i+1}S_n)H_{k+1}(\Delta^{i-1}S_{n+1})}{H_{k+1}(\Delta^iS_{n+1})
H_{k+1}(\Delta^iS_{n})}.\label{h3}
\end{eqnarray}

Then, from the formulae \eqref{h1}--\eqref{h3}, we have
\begin{eqnarray}
\prod_{i=1}^{m}(\varepsilon_{(m+1)k+i,m}^{(n+1)}-\varepsilon_{(m+1)k+i,m}^{(n)})
&\!\!=\!\!&-\frac{H_{k+1}(\Delta^2S_{n})H_{k}(\Delta^{m+1}S_{n+1})}
{H_{k+1}(\Delta S_{n+1})H_{k+1}(\Delta S_{n})}\prod_{i=2}^m \frac{H_{k+1}(\Delta^{i+1}S_{n})
H_{k+1}(\Delta^{i-1}S_{n+1})}
{H_{k+1}(\Delta^iS_{n})H_{k+1}(\Delta^iS_{n+1})}\nonumber\\
&\!\!=\!\!&-\frac{H_{k}(\Delta^{m+1}S_{n+1})H_{k+1}(\Delta^{m+1}S_{n})}
{H_{k+1}(\Delta S_{n})H_{k+1}(\Delta^mS_{n+1})}.\label{r1}
\end{eqnarray}
Comparing (\ref{l1}) and (\ref{r1}), we obtain the rule (\ref{gd}) of the multistep $\varepsilon$--algorithm where
the lower index $k$ is replaced by $(m+1)k+1$.

We can also derive the following formula
\begin{eqnarray}
\prod_{i=2}^{m+1}(\varepsilon_{(m+1)k+i,m}^{(n+1)}-\varepsilon_{(m+1)k+i,m}^{(n)})
&\!\!=\!\!&\frac{H_{k+2}(\Delta S_{n})H_{k+1}(\Delta^{m}S_{n+1})}
{H_{k+1}(\Delta^{m+1}S_{n+1})H_{k+1}(\Delta^{m+1}S_{n})}\prod_{i=2}^m \frac{H_{k+1}(\Delta^{i+1}S_{n})
H_{k+1}(\Delta^{i-1}S_{n+1})}
{H_{k+1}(\Delta^iS_{n})H_{k+1}(\Delta^iS_{n+1})}\nonumber\\
&\!\!=\!\!&\frac{H_{k+2}(\Delta S_{n})H_{k+1}(\Delta S_{n+1})}
{H_{k+1}(\Delta^2 S_{n})H_{k+1}(\Delta^{m+1}S_{n+1})}.\label{r2}
\end{eqnarray}
Comparing (\ref{l2}) and (\ref{r2}), we obtain the rule (\ref{gd}) of the multistep $\varepsilon$--algorithm where the lower
index $k$ is replaced by $(m+1)k+2$.

Besides, we have, for $i=2,\ldots,m$,
\begin{eqnarray}
&&\prod_{j=i+1}^{m+i}(\varepsilon_{(m+1)k+j,m}^{(n+1)}-\varepsilon_{(m+1)k+j,m}^{(n)}) =  \nonumber\\
&&\mbox{~~~~~~~~~~~~~~~~~~~~~} \prod_{j=i+1}^{m}(\varepsilon_{(m+1)k+j,m}^{(n+1)}-\varepsilon_{(m+1)k+j,m}^{(n)})
\prod_{j=2}^{i-1}(\varepsilon_{(m+1)(k+1)+j,m}^{(n+1)}-\varepsilon_{(m+1)(k+1)+j,m}^{(n)})\nonumber\\
&&\mbox{~~~~~~~~~~~~~~~~~~~~~~~~~~~~~~~}
\cdot(\varepsilon_{(m+1)(k+1),m}^{(n+1)}-\varepsilon_{(m+1)(k+1),m}^{(n)})
(\varepsilon_{(m+1)(k+1)+1,m}^{(n+1)}-\varepsilon_{(m+1)(k+1)+1,m}^{(n)}) \nonumber\\
&&\mbox{~~~~~~~~~~~~~~~~~~~~~} = -\prod_{j=i+1}^m\frac{H_{k+1}(\Delta^{j+1}S_{n})H_{k+1}(\Delta^{j-1}S_{n+1})}
{H_{k+1}(\Delta^j S_{n+1})H_{k+1}(\Delta^j S_{n})}\prod_{j=2}^{i-1}
\frac{H_{k+2}(\Delta^{j+1}S_{n})H_{k+2}(\Delta^{j-1}S_{n+1})}
{H_{k+2}(\Delta^jS_{n+1})H_{k+2}(\Delta^jS_{n})}\nonumber\\
&&\mbox{~~~~~~~~~~~~~~~~~~~~~~~~~~~~~~~}
 \cdot\frac{H_{k+2}(\Delta S_{n})H_{k+1}(\Delta^{m}S_{n+1})}
{H_{k+1}(\Delta^{m+1}S_{n+1})H_{k+1}(\Delta^{m+1}S_{n})}\cdot\frac{H_{k+2}(\Delta^2S_{n})
H_{k+1}(\Delta^{m+1}S_{n+1})}
{H_{k+2}(\Delta S_{n+1})H_{k+2}(\Delta S_{n})} \nonumber\\
&&\mbox{~~~~~~~~~~~~~~~~~~~~~}
=-\frac{H_{k+1}(\Delta^i S_{n+1})H_{k+2}(\Delta^{i}S_{n})}
{H_{k+1}(\Delta^{i+1}S_{n})H_{k+2}(\Delta^{i-1}S_{n+1})}.\label{r3}
\end{eqnarray}
Comparing (\ref{l3}) and (\ref{r3}), we obtain the rule (\ref{gd}) of the multistep $\varepsilon$--algorithm with the lower
index $k$ replaced by $(m+1)k+i+1$, for $ i=2, 3, \ldots,m$.
Therefore, the multistep $\varepsilon$--algorithm has been derived from the definition (\ref{slt1}), (\ref{slt2}) and (\ref{slt3})
of the multistep Shanks' transformation, and we have the

\begin{theorem}
~~\\
The multistep Shanks' transformation defined by (\ref{slt1}), (\ref{slt2}), and (\ref{slt3}) can be implemented by the recursive rules \eqref{gd} of the multistep $\varepsilon$--algorithm, with the initializations \eqref{ini}.
\end{theorem}

\section{An extended discrete Lotka--Volterra system}
\label{slv}

Recently, as explained in Section \ref{sce}, it has been shown that integrable systems are closely
related to numerical algorithms. On one hand, some numerical
algorithms are found to be soliton equations. For example, one step
of the $QR$--algorithm is equivalent to the time evolution of the
finite non--periodic Toda lattice \cite{1}. The $\varepsilon$--algorithm
is nothing but the fully--discrete potential KdV equation, and the
$\rho$--algorithm is considered to be the fully--discrete cylindrical
KdV equations or the Milne--Thomson equation, see \cite{wynn,int4,naga, 5,vein}. On the other hand,
integrable systems can be used for designing new numerical
algorithms. For example, the discrete Lotka--Volterra system has
applications in numerical algorithms for computing singular values
\cite{6, 7, 8}, the continuous--time Toda equation leads to a new
algorithm for computing the Laplace transform of a given analytic
function \cite{9}, and the discrete relativistic Toda molecule equation leads to a new Pad\'{e}
approximation algorithm for formal power series \cite{mine}.

\vskip 2mm

In this section, we will show that there exist a {\it Miura transformation} between the multistep $\varepsilon$--algorithm
\eqref{gd} and a discrete integrable system. In fact, if we set $\left(a^{(n)}_{k-\frac{m-1}{2}}\right)^{-1}=
\varepsilon_{k,m}^{(n+1)}-\varepsilon_{k,m}^{(n)}$, then equation \eqref{gd} is transformed into
the {\it extended discrete Lotka--Volterra equation}
\begin{equation}
\label{exlv1}
\prod_{i=0}^{m-1}a_{k-\frac{m-1}{2}+i}^{(n+1)}-\prod_{i=0}^{m-1}a_{k-\frac{m-1}{2}+i}^{(n)}=\frac{1}
{a_{k+\frac{m+1}{2}}^{(n)}}-\frac{1}{a_{k-\frac{m+1}{2}}^{(n+1)}}.
\end{equation}
This equation can be considered as the time discretization, for $N=-1$, of
\begin{equation}
\frac{d}{dt}\left(\prod_{i=0}^{m-1}a_{k-\frac{m-1}{2}+i}\right)=\prod_{i=0}^{-N-1}a^{-1}_{k+\frac{m+1}{2}+i}-
\prod_{i=0}^{-N-1}a^{-1}_{k-\frac{m+1}{2}-i},\quad
m=1,2,\ldots, \quad N=-1,-2,\ldots,\label{exlv2}
\end{equation}
which is called the {\it extended Lotka--Volterra equation}. This equation was first proposed in \cite{Nar},
and it was developed in \cite{Hu}. Indeed, with $N=-1$, \eqref{exlv2} becomes
\begin{equation}
\frac{d}{dt}\left(\prod_{i=0}^{m-1}a_{k-\frac{m-1}{2}+i}\right)=\frac{1}{a_{k+\frac{m+1}{2}}}-
\frac{1}{a_{k-\frac{m+1}{2}}}. \label{g1}
\end{equation}
Now, consider $n$ as the discretization of $t$, and replace the derivative in the left hand side of \eqref{g1}
by the forward difference $\Delta$ acting on $n$. The left hand side becomes
$$\prod_{i=0}^{m-1}a_{k-\frac{m-1}{2}+i}^{(n+1)}-\prod_{i=0}^{m-1}a_{k-\frac{m-1}{2}+i}^{(n)}.$$
Then, replace $a_k$ in the first term of the right hand side of \eqref{g1} by $a_k^{(n)}$,and, in its second term,
by $a_k^{(n+1)}$. We get \eqref{exlv1}.

\vskip 2mm

Using the relations \eqref{bt2} and \eqref{bt1}, we obtain the solution of \eqref{exlv1}
\begin{eqnarray*}
a_{(m+1)k+1-\frac{m-1}{2}}^{(n)}=-\frac{1}{\varepsilon_{(m+1)k+1,m}^{(n+1)}-\varepsilon_{(m+1)k+1,m}^{(n)}}=
-\frac{F^n_{(m+1)k+1}F^{n+1}_{(m+1)k+1}}{F^n_{(m+1)k+2}F^{n+1}_{(m+1)k}},&&\\
a_{(m+1)k+i-\frac{m-1}{2}}^{(n)}=\frac{1}{\varepsilon_{(m+1)k+i,m}^{(n+1)}-\varepsilon_{(m+1)k+i,m}^{(n)}}=
\frac{F^n_{(m+1)k+i}F^{n+1}_{(m+1)k+i}}{F^n_{(m+1)k+i+1}F^{n+1}_{(m+1)k+i-1}},&&\\
i=2,\ldots,m+1,&&
\end{eqnarray*}
that is
\begin{eqnarray}
a_{(m+1)k-\frac{m-1}{2}}^{(n)}&=&\frac{H_{k}(\Delta^{m+1}S_n)H_{k}(\Delta^{m+1}S_{n+1})}{H_{k+1}(\Delta S_n)
H_{k}(\Delta^{m}S_{n+1})},\label{exlvsol1}\\
a_{(m+1)k-\frac{m-1}{2}+1}^{(n)}&=&-\frac{H_{k+1}(\Delta S_n)H_{k+1}(\Delta S_{n+1})}{H_{k+1}(\Delta^2 S_n)
H_{k}(\Delta^{m+1}S_{n+1})},\label{exlvsol2}\\
a_{(m+1)k-\frac{m-1}{2}+j}^{(n)}&=&\frac{H_{k+1}(\Delta^j S_n)H_{k+1}(\Delta^j S_{n+1})}{H_{k+1}(\Delta^{j+1}S_n)
H_{k+1}(\Delta^{j-1}S_{n+1})},\label{exlvsol3}
\end{eqnarray}
where $j=2,\ldots,m$, and $k=-m+1,-m+2,\ldots$, with the initial values
\begin{eqnarray}
a^{(n)}_{-m-\frac{m-1}{2}}=\infty, a^{(n)}_{-m+1-\frac{m-1}{2}}=\cdots=a^{(n)}_{-1-\frac{m-1}{2}}=n, a^{(n)}_{-\frac{m-1}{2}}
=1/\Delta S_n.\label{exlvini}
\end{eqnarray}

The difference equation \eqref{exlv1}, with the initial values \eqref{exlvini}, is said to be the integrable time discretization
of the extended Lotka--Volterra equation \eqref{exlv2} in the sense that its solution is given by
\eqref{exlvsol1}--\eqref{exlvsol3}.
Conversely,  the extended discrete Lotka--Volterra equation \eqref{exlv2}
can be seen as the time continuation of \eqref{exlv1} with the initializations \eqref{exlvini}.

\vskip 2mm

Consider the particular case $m=1$. Then, \eqref{gd} reduces to the $\varepsilon$--algorithm and equation
\eqref{exlv1} becomes
\begin{eqnarray*}
a_{k}^{(n+1)}-a_k^{(n)}&=&\frac{1}{a_{k+1}^{(n)}}-\frac{1}{a_{k-1}^{(n+1)}}.
\end{eqnarray*}

By the dependent variable transformation
\begin{eqnarray*}
\frac{u_k^{(n)}}{u_{k-1}^{(n+1)}}=\frac{a_{k-1}^{(n+1)}}{a_{k+1}^{(n)}},
\end{eqnarray*}
we obtain the discrete Lotka--Volterra equation
\begin{equation}
\label{lv}
u_{k}^{(n+1)}\left(1+u_{k-1}^{(n+1)}\right)=u_k^{(n)}\left(1+u_{k+1}^{(n)}\right).
\end{equation}
Then, the $\varepsilon$--algorithm can be transformed into the discrete Lotka--Volterra equation \eqref{lv} through
the following Miura transformation
\begin{equation}
\nonumber
\frac{u_k^{(n)}}{u_{k-1}^{(n+1)}}=\frac{\varepsilon_{k+1,1}^{(n+1)}-\varepsilon_{k+1,1}^{(n)}}
{\varepsilon_{k-1,1}^{(n+2)}-\varepsilon_{k-1,1}^{(n+1)}}.
\end{equation}
Thus, the $\varepsilon$--algorithm can be considered as the discrete Lotka--Volterra equation \eqref{lv}, and
more generally, the multistep $\varepsilon$--algorithm \eqref{gd} is equivalent to the extended discrete
Lotka--Volterra equation \eqref{exlv1}.

\section{Conclusion and future researches}
\label{ssfutur}

Starting from the recursive rule (\ref{gd}) of the multistep $\varepsilon$--algorithm, we first obtained,
from the Hirota's bilinear method, the coupled relations (\ref{bt2})--(\ref{bt3}). Then, applying the Sylvester's identity to
the determinants
${H_{k}(\Delta^iS_n)}$ and ${\Phi_{k}(\Delta^iS_n)}$, we got the formulae \eqref{bf1}--\eqref{bf2}
which express the quantities $\varepsilon_{k,m}^{(n)}$ as ratios of determinants.
Thus, we were able to {\it prove} that
the $\varepsilon_{k,m}^{(n)}$'s are defined as ratios of determinants, and then to
{\it derive} the recursive rule (\ref{gd}) of the multistep $\varepsilon$--algorithm, with the initializations \eqref{ini}, from the determinantal formulae defining
the quantities $\varepsilon_{k,m}^{(n)}$. It must be noticed that, contrarily to the approaches of \cite{wynn} and
\cite{huweni}, we did not make use of the Schweins' determinantal identity, but only of the Sylvester's one.
The difficult point was to find to which determinants this identity had to be applied.
Then, we showed that the multistep $\varepsilon$--algorithm was related to an extended discrete Lotka--Volterra system.

When $m=1$, the relations \eqref{slt3} disappear, and the Hirota's bilinear method leads to a new proof
that the $\varepsilon$--algorithm
of Wynn implements the Shanks' sequence transformation and, reciprocally, that the
quantities computed by this algorithm are expressed by the ratios of Hankel determinants defining the Shanks'
transformation.

\vskip 2mm

The approach developed above could possibly be extended to other nonlinear convergence acceleration
algorithms such as, for example, the $q$-difference version of the $\varepsilon$--algorithm proposed in \cite{he},
or its two generalizations given in \cite{geneps}, or the other one presented in \cite{hh},  or the general
$\varepsilon$--algorithm of \cite{5cars}, or the $\rho$--algorithm  \cite{5rho}, and the $\gamma$--algorithm which
generalizes it \cite{th}.
Other algorithms related to them, such as the $qd$, the $\eta$, the $\omega$, and the $rs$--algorithms, and
the $g$--decomposition, could also possibly be treated in a similar way (see \cite{cbmrz} for their definitions).
The quantities computed by these algorithms are all defined as ratios of determinants.
These extensions, as well as extensions to other acceleration algorithms, will be the subject of future works.
Let us mention that the confluent form of the multistep $\varepsilon$--algorithm is studied in \cite{mepsc}.
It leads to a multistep Lotka--Volterra equation.

\vskip 2mm

\noindent{\bf Acknowledgements:}
This work was partially
supported by the National Natural Science Foundation of China (Grant
no. 11071241), and the knowledge innovation program of LSEC and the
Institute of Computational Math., AMSS, CAS. C. Brezinski would like to thanks X.B. Hu, the State Key Laboratory
of Scientific and Engineering
Computing (LSEC), and the Institute of Computational Mathematics, AMSS, CAS, for inviting him for a stay during which
part of this work was done.
The work of Michela Redivo--Zaglia was partially supported by MIUR, PRIN grant no. 20083KLJEZ-003, and by University of Padova, Project 2008 no. CPDA089040.

\end{document}